# Adaptive FEM-Based Nonrigid Image Registration Using Truncated Hierarchical B-splines

Aishwarya Pawar[a], Yongjie Zhang[a,*], Yue Jia[b], Xiaodong Wei[a], Timon Rabczuk[b], Chiu Ling Chan[b], Cosmin Anitescu[b]

[a]*Department of Mechanical Engineering, Carnegie Mellon University, USA*
[b]*Institute of Structural Mechanics, Bauhaus-Universität Weimar, Germany*

**Abstract**

We present an efficient approach of Finite Element Method (FEM)-based nonrigid image registration, in which the spatial transformation is constructed using truncated hierarchical B-splines (THB-splines). The image registration framework minimizes an energy functional using an FEM-based method and thus involves solving a large system of linear equations. This framework is carried out on a set of successively refined grids. However, due to the increased number of control points during subdivision, large linear systems are generated which are generally demanding to solve. Instead of using uniform subdivision, an adaptive local refinement scheme is carried out, only refining the areas of large change in deformation of the image. By incorporating the key advantages of THB-spline basis functions such as linear independence, partition of unity and reduced overlap into the FEM-based framework, we improve the matrix sparsity and computational efficiency. The performance of the proposed method is demonstrated on 2D synthetic and medical images.

*Keywords:* Nonrigid Image Registration, Truncated Hierarchical B-splines, Adaptive Local Refinement, Finite Element Method

## 1. Introduction

Image registration has significant applications in the field of remote sensing, computer vision, medical image analysis and other industrial applications [1, 2]. The technique basically involves aligning two or more images and finding the accurate correspondence between them.

*Corresponding Author
Email address:* `jessicaz@andrew.cmu.edu` (Yongjie Zhang)

This plays a critical role in integrating useful data, detecting changes in certain important features and thus enabling a better understanding of the image data. A comprehensive and exhaustive survey on the various image registration methods, their classification and application-specific advantages can be found in [3, 4].

We can fundamentally break down the entire image registration process into the following salient steps. The first step involves computing the similarity between two input images. The given source image undergoes a deformation in order to match with the target image. We compute a metric that measures the similarity between them and drives the registration to maximize the similarity. A spatial transformation function is defined in a parametric form and the optimal parameters are computed by maximizing the similarity between the source and target images.

Based on the type of image deformation, registration methods are classified as rigid and nonrigid. Under rigid registration, simple transformations such as translation, rotation and reflection are performed. In nonrigid registration it is possible to model more complex deformations. However there is a need to impose certain regularization constraints to control the smoothness of the deformation and the convergence of the solution. In spline-based registration [5, 6, 7], a control grid overlays the image and an optimum spatial transformation is constructed. As compared to the standard isoparametric basis functions, representing the deformation and computing the spatial transformation with hierarchical B-splines (HB-splines) or truncated hierarchical B-splines (THB-splines) is more accurate for the same number of degrees of freedom. To capture large scale deformations, coarser grids are used with large support. To capture highly localized deformations, finer grids are used. The properties of B-splines such as local control, smoothness and compact support are used to model the deformation in the image.

FEM-based image registration is one of the most promising image registration methods. In FEM-based methods, the numerical implementation is carried out by solving an optimization problem using a system of linear equations. In [8, 9], an energy functional is constructed based on the sum of squared differences in the intensity values, and then solved using $L^2$-gradient flow method. Regularization constraints are also added to ensure that the deformation is smooth. FEM-based image registration methods involve the assembly of large matrices and increased computational complexity. Furthermore, the implementation is carried out on a hierarchy of uniform grids, in order to capture coarse deformations on coarser grids and fine-scale deformations on finer grids. This makes the assembly of the matrices more complex and cumbersome.

To resolve these matrix assembly issues, we can reduce the size of the matrices and increase the sparsity of the matrices. In this paper, we introduce adaptive local refinement of the control grid instead of carrying out uniform subdivision, enabling highly localized fine-scaled deformations only in desired regions of the image domain. Local refinement can be achieved using HB-splines, which have been previously used in [10, 11] to solve image registration problems. Pertaining to FEM-based image registration problems, we employ THB-splines for adaptive local refinement. In addition to being linearly independent and non-negative, these basis functions have smaller support and form a partition of unity. These properties effectively improve the computational efficiency in solving FEM-based problems and help in achieving the desired accuracy using fewer control points by local refinement and sparser matrices.

The rest of the paper is organized as follows. In section 2, we first present the framework of image registration and introduce our adaptive grid refinement strategy using THB-splines in section 3. In section 4, we apply our approach to solve numerical examples in the form of 2D synthetic and medical images, by comparing the results in terms of computational time and accuracy with the uniform grid refinement strategy. The paper ends with a brief summary in section 5.

## 2. Image Registration Model

In this section we briefly review the registration framework, in which the formulation of the energy functional and FEM-based method is explained. We then review THB-splines and introduce our adaptive local refinement scheme, describing our motivation of using THB-splines by explaining their key advantages over using HB-spline basis functions.

### 2.1. Registration Framework

The process of image registration involves computation of a spatial transformation function $f(\mathbf{x})$, that aligns the source image $I_1(\mathbf{x})$ to the target image $I_2(\mathbf{x})$. Bi-cubic B-splines are used to define the spatial transformation and create a $C^2$-continuous mapping from $R^2 \to R^2$. The spatial transformation function $f(\mathbf{x})$ is defined as:

$$f(\mathbf{x}) = \sum_{k=1}^{N_b} P_k \phi_k(\mathbf{x}), \tag{1}$$

where $P_k$ is a set of control points associated with the bivariate basis functions $\phi_k(\mathbf{x})$. $N_b$ represents the total number of basis functions. $\phi_k(\mathbf{x})$ is the tensor product of univariate B-spline basis functions $N_{i,p}(\xi)$ and $N_{j,q}(\eta)$ defined on the knot vectors $\Xi = \{\xi_1, \cdots, \xi_{n_1+p+1}\}$ and $\Upsilon = \{\eta_1, \cdots, \eta_{n_2+q+1}\}$ in $\xi$ and $\eta$ directions, respectively. $n_1$ and $n_2$ are the number of univariate basis functions. $p$ and $q$ are the degree of polynomials in $\xi$ and $\eta$ directions, respectively. This is written as

$$\phi_k(\mathbf{x}) = N_{i,p}(\xi) N_{j,q}(\eta), \tag{2}$$

where $i = 1, 2, ..., n_1$ and $j = 1, 2, ..., n_2$. The detailed description on B-splines along with their applications in the field of computational modeling and isogeometric analysis can be found in [12, 13]. At the start of the registration process we define the initial location of the control points so as to create an identity map, $f(\mathbf{x}) = \mathbf{x}$. This also represents the source image $I_1(\mathbf{x})$. Then, we proceed to find the optimal transformation $f(\mathbf{x})$ such that $I_1(f(\mathbf{x})) \approx I_2(\mathbf{x})$. The transformation of the source image is driven by minimizing an energy functional, which uses the sum of squared differences ($SSD$) of the two images as the driving force to align them. The detailed description of the formulation of the energy functional can

2.2 Registration Framework

be found in [8, 9]. The energy functional is given as:

$$E(f(\mathbf{x})) = \int_\Omega g(\mathbf{x})(I_2(\mathbf{x}) - I_1(f(\mathbf{x})))^2 \, d\Omega + \lambda_1 \int_\Omega (\|f_{,\xi}(\mathbf{x})\|_2^2 + \|f_{,\eta}(\mathbf{x})\|_2^2) \, d\Omega \\ + \lambda_2 \int_\Omega (\|f_{,\xi}(\mathbf{x})\|_2^2 \|f_{,\eta}(\mathbf{x})\|_2^2 - (<f_{,\xi}(\mathbf{x}), f_{,\eta}(\mathbf{x})>)^2) d\Omega, \quad (3)$$

where $\lambda_1$ and $\lambda_2$ are the regularization parameters set to a constant value during the entire registration process. $f_{,\xi}(\mathbf{x})$ and $f_{,\eta}(\mathbf{x})$ are the first derivatives of $f(\mathbf{x})$ in $\xi$ and $\eta$ directions respectively and $<f_{,\xi}(\mathbf{x}), f_{,\eta}(\mathbf{x})>$ is the inner product operator which equals to $f_{,\xi}(\mathbf{x})^T f_{,\eta}(\mathbf{x})$. $g(\mathbf{x})$ is defined as

$$g(\mathbf{x}) = \frac{1}{\sqrt{\gamma + (\nabla_\xi I_1(f(\mathbf{x})))^2 + (\nabla_\eta I_1(f(\mathbf{x})))^2}}, \quad (4)$$

as given in [8]. $g(\mathbf{x})$ is used to accelerate the registration process in the homogeneous regions and slow down the registration in the inhomogeneous regions of the image. $\gamma$ is a small number introduced to prevent the division by zero. Here we set $\gamma$ as $10^{-12}$.

The first term in Equation (3) drives the registration process by minimizing the differences between the images. Unlike rigid deformations, non-rigid deformations can achieve free-form deformations of larger magnitude. But this sometimes results in unrealistic changes within the image. Regularization is therefore very essential in ensuring the smooth deformation of the image. The second and third terms in Equation (3), called the regularization constraints, are used to prevent these unrealistic changes and also make the transformation as smooth as possible. The first regularization term ensures smooth variation of $f(\mathbf{x})$ in $\xi$ and $\eta$ directions, and thus prevents large deformation of the evolving image. The second regularization term ensures the consistency of the area element during deformation. This ensures that the deformation does not result in unrealistic results and converges to the desired accuracy.

As given in [9], the $L^2$-gradient flow method is chosen to convert the optimization problem to an initial value problem of an ordinary differential equation. By minimizing the energy functional with respect to the spatial transformation function, we compute the most optimal $f(\mathbf{x})$ that maximizes the match between the images. Since the energy functional is a function of $f(\mathbf{x})$, we differentiate $E(\mathbf{x})$ with respect to $f(\mathbf{x})$ evaluated at the previous time step and update it according to the following equation,

$$\frac{df(\mathbf{x})}{dt} = -\delta E(f(\mathbf{x})). \quad (5)$$

To get the weak formulation, we use a test function $\Psi(\mathbf{x})$ and integrate over the image domain. $\Psi(\mathbf{x})$ is defined to be the equal to $\phi(\mathbf{x})$ according to Bubnov-Galerkin method. The weak formulation is given as

$$\int_\Omega \frac{df(\mathbf{x})}{dt} \Psi(\mathbf{x}) d\Omega = -\delta E(f(\mathbf{x}), \Psi(\mathbf{x})), \quad (6)$$

where $\delta E(f(\mathbf{x}), \Psi(\mathbf{x}))$ is the first-order variation of the energy functional defined in Equation (3). Thus we can write

$$\begin{aligned}\delta E(f(\mathbf{x}), \Psi(\mathbf{x})) &= -2\int_{\Omega} g(\mathbf{x})(I_2(\mathbf{x}) - I_1(f(\mathbf{x})))\nabla I_1(f(\mathbf{x}))\Psi(\mathbf{x})\mathrm{d}\Omega \\ &\quad +\lambda_1 \int_{\Omega} 2(f_{,\xi}(\mathbf{x})\Psi_{,\xi}(\mathbf{x}) + f_{,\eta}(\mathbf{x})\Psi_{,\eta}(\mathbf{x}))\mathrm{d}\Omega \\ &\quad +\lambda_2 \int_{\Omega} 2(f_{,\xi}(\mathbf{x})\|f_{,\eta}(\mathbf{x})\|^2 \Psi_{,\xi}(\mathbf{x}) + \|f_{,\xi}(\mathbf{x})\|^2 f_{,\eta}(\mathbf{x})\Psi_{,\eta}(\mathbf{x}) \\ &\quad -\langle f_{,\xi}(\mathbf{x}), f_{,\eta}(\mathbf{x})\rangle f_{,\eta}(\mathbf{x})\Psi_{,\xi}(\mathbf{x}) \\ &\quad -f_{,\xi}(\mathbf{x})\langle f_{,\eta}(\mathbf{x}), f_{,\xi}(\mathbf{x})\rangle \Psi_{,\eta}(\mathbf{x}))\mathrm{d}\Omega,\end{aligned} \quad (7)$$

where the derivation of the first-order variation of the energy functional can be found in [9].

Since $f(\mathbf{x})$ can be written as shown in Equation (1), we can write

$$\sum_{i=1}^{N_b} \int_{\Omega} \frac{P_i^{t+1} - P_i^t}{\epsilon} \phi_i(\mathbf{x})\Psi_j(\mathbf{x})\mathrm{d}\Omega = -\delta E_j(f(\mathbf{x}), \Psi(\mathbf{x})), \quad (8)$$

where $j = 1, 2, ..., N_b$, $N_b$ representing the total number of bivariate B-splines. $\epsilon$ denotes the time step. In matrix form we rewrite the above equation as

$$\mathbf{M}(\mathbf{P}^{t+1} - \mathbf{P}^t) = -\epsilon \mathbf{E}, \quad (9)$$

where in matrix $\mathbf{M}$, $M_{ji} = \int_{\Omega} \phi_i \Psi_j \mathrm{d}\Omega$ ($i, j = 1, \cdots, N_b$). $\mathbf{P}^t$ and $\mathbf{P}^{t+1}$ represent the matrices of the control point coordinates at the previous and the current time step respectively. $\mathbf{E}$ is the vector storing $\delta E_j(f(\mathbf{x}), \Psi(\mathbf{x}))$ where

$$\begin{aligned}\delta E_j(f(\mathbf{x}), \Psi(\mathbf{x})) &= -2\int_{\Omega} g(\mathbf{x})(I_2(\mathbf{x}) - I_1(f(\mathbf{x})))\nabla I_1(f(\mathbf{x}))\Psi_j(\mathbf{x})\mathrm{d}\Omega \\ &\quad +\lambda_1 \int_{\Omega} 2(f_{,\xi}(\mathbf{x})\Psi_{j,\xi}(\mathbf{x}) + f_{,\eta}(\mathbf{x})\Psi_{j,\eta}(\mathbf{x}))\mathrm{d}\Omega \\ &\quad +\lambda_2 \int_{\Omega} 2(f_{,\xi}(\mathbf{x})\|f_{,\eta}(\mathbf{x})\|^2 \Psi_{j,\xi}(\mathbf{x}) + \|f_{,\xi}(\mathbf{x})\|^2 f_{,\eta}(\mathbf{x})\Psi_{j,\eta}(\mathbf{x}) \\ &\quad -\langle f_{,\xi}(\mathbf{x}), f_{,\eta}(\mathbf{x})\rangle f_{,\eta}(\mathbf{x})\Psi_{j,\xi}(\mathbf{x}) \\ &\quad -f_{,\xi}(\mathbf{x})\langle f_{,\eta}(\mathbf{x}), f_{,\xi}(\mathbf{x})\rangle \Psi_{j,\eta}(\mathbf{x}))\mathrm{d}\Omega.\end{aligned} \quad (10)$$

Equation (9) is solved using an FEM-based method. Numerical integration is carried out by Gaussian quadrature rule of order 6. We can see that the fast computation of $\mathbf{M}$ is crucial for achieving accurate results efficiently.

Similarity metrics, such as the mean squared difference ($MSD$) and the similarity ratio ($RS$) [8], provide a way to quantitatively evaluate the results of the registration process. $RS$ is defined as

$$RS(I_2(\mathbf{x}), I_1^t(f(\mathbf{x}))) = \left(1 - \frac{\|I_2(\mathbf{x}) - I_1^t(f(\mathbf{x}))\|_{L^2}}{\|I_2(\mathbf{x}) - I_1(\mathbf{x})\|_{L^2}}\right) \times 100, \quad (11)$$

where $I_1^t(f(\mathbf{x}))$ is the image obtained at the $t^{th}$ time step. When $RS = 100\%$, this corresponds to a perfect registration.

2.2. Truncated Hierarchical B-splines (THB-splines)

THB-splines are developed based on HB-spline basis functions. Let us first briefly review HB-spline basis functions. We consider univariate basis functions defined on the open knot vector $\Xi = \{\xi_1, \cdots, \xi_{n_1+p+1}\}$, where $n_1$ and $p$ are the number of univariate basis functions and the degree of the polynomials in the $\xi$ direction, respectively. The local support of the basis function $N_{i,p}(\xi)$ is $[\xi_i, \xi_{i+p+1}]$, denoted by $supp(N_{i,p}(\xi))$. The maximum refinement level is denoted as $l_{max}$. The knot vectors defined at a particular refinement level are obtained by bisecting the knot vectors of the previous refinement level. A basis function $N_{i,p}^l(\xi)$ at a given refinement level $l$ can be represented as a linear combination of a subset of the basis functions from the next refinement level $(l+1)$, which are the children basis functions of $N_{i,p}^l(\xi)$. We have

$$N_{i,p}^l(\xi) = \sum_{k=0}^{N_c-1} S_{k,p} N_{k,p}^{l+1}(\xi), \quad (12)$$

where $N_{k,p}^{l+1}(\xi)$ are the children basis functions and $N_c$ is the number of $N_{k,p}^{l+1}(\xi)$. We obtain the refinement coefficients $S_{k,p}$ using the Oslo Algorithm [14]. In the similar manner, bivariate basis functions $\phi_k^l(\mathbf{x})$ (Equation (2)) can be represented as a linear combination of their children basis functions, which are completely contained in the local support of $\phi_k^l(\mathbf{x})$, $[\xi_i, \xi_{i+p+1}] \times [\eta_j, \eta_{j+q+1}]$. Details of local refinement using HB-splines can be found in [15], in which efficient ways to construct analysis-suitable HB-splines is described.

HB-splines possess desirable properties such as linear independence, non-negativity and support for local refinement. However, the basis functions need to be rationalized to form a partition of unity, which is complicated and time-consuming. Moreover, HB-splines have extensive overlapping among coarser and finer splines. This overlap increases with increasing levels of refinement. To develop polynomial basis functions satisfying partition of unity and reduce the overlapping of the B-splines from different levels, THB-splines were proposed [16]. We define the total support of refined B-splines basis functions $\phi_r^l(\mathbf{x})$ at a particular refinement level $l$ as $\Omega^{l+1} = \cup supp(\phi_r^l(\mathbf{x}))$. Equation (12) is modified for THB-splines as follows

$$trunc(\phi_k^l(\mathbf{x})) = \sum_{supp(\phi_p^{l+1}(\mathbf{x})) \not\subset \Omega^{l+1}} S_{k,p} \phi_p^{l+1}(\mathbf{x}), \quad (13)$$

where the truncated B-splines at a particular refinement level are represented using their children basis functions that are not present in the hierarchical B-spline basis. This enables local refinement of the selected basis functions at different refinement levels $(l = 1, \cdots, l_{max})$, further reducing the overlapping between the basis functions and satisfying partition of unity. We explain local refinement based on THB-spline basis functions on two consecutive levels, $l$ and $l+1$, in the following steps:

1. At a particular refinement level $l$, a subset of B-spline basis functions satisfying the refinement criterion are identified, which form the set $\phi_r^l$. The remaining B-spline basis functions are set as active ($\phi_a^l$). Among the active basis functions, the truncated basis functions are represented using Equation (13).
2. The children B-spline basis functions of all the basis functions $\phi_r^l$ are set as active ($\phi_a^{l+1}$).
3. The process ends by collecting all the active B-spline basis functions at levels $l$ and $l+1$, in order to get the THB-splines basis functions, $\phi_{thb}^{l+1}$.

Thus we can obtain
$$\phi_{thb}^{l+1} = \phi_a^l \cup \phi_a^{l+1}. \tag{14}$$

The refinement procedure is implemented in a recursive manner until we reach the refinement level $l_{max}$.

## 3. Adaptive Local Refinement

In this section, we describe the procedure to implement adaptive local refinement using THB-splines. A multilevel technique was proposed in [8], where the entire registration process is carried out on a set of grids that are uniformly subdivided for each refinement level. It was observed that the similarity metric does not continue to increase further or even drops after certain number of iterations. This happens when we compute the spatial transformation on a single grid, resulting in distortion of certain elements in the grid due to large deformations. Instead of uniformly refining the grid for each refinement level, we can use the information from the evolving and target images, and refine only those regions that undergo large change in image deformation. Starting from a particular refinement level, we carry out adaptive local refinement and reinitialize the control mesh. The image obtained at the previous refinement level is used as a source image for the registration at the next refinement level. This process continues until the maximum refinement level is reached. We use the magnitude of the gradient of the difference in the two images $I_1(\mathbf{x})$ and $I_2(\mathbf{x})$, $I_g = |\nabla(I_1(f(\mathbf{x})) - I_2(\mathbf{x}))|$, to compute the refinement. We can detect the areas of large change in deformation in this manner and locally refine them. The size of the elements at the first level is decided based on the complexity of the image. For medical images, the initial mesh size should not be too coarse as this may affect the convergence of the registration process. Even if the initial registration is faster, by using coarser meshes it can sometimes cause unrealistic deformations within the image. Thus to maintain the stability of the entire registration process, the initial size of the elements should be kept sufficiently fine in order to achieve accurate registration results.

The detailed algorithm is described in the following: Given a pair of images, $I_1(\mathbf{x})$ and $I_2(\mathbf{x})$. The process starts at level-1, which is a uniform grid. The initial set of control points $P_{thb}^1$ are defined such that $f(\mathbf{x}) = \sum_{i=1}^{N} P_{thb,i}^1 \phi_{thb,i}^1 = \mathbf{x}$. $N$ is the total number of basis functions at level-1. For $l = 1, ..., l_{max}$,

1. Compute $I_g = |\nabla(I_1^{l-1}(f(\mathbf{x})) - I_2(\mathbf{x}))|$ at each control grid of level-$l$ if $l > 1$. If $l = 1$, then we directly proceed to step 3. Here $I_1^{l-1}(f(\mathbf{x}))$ is the evolving image obtained from the previous level-$(l-1)$. Let $G_{mean}$ be the average value of $I_g$.
2. Loop over the active B-spline basis functions $\phi_j^{l-1}$, $j = 1, 2, ..., N_a$, where $N_a$ is the total number of active B-spline basis functions at level-$(l-1)$. Loop over $supp(\phi_j^{l-1})$, and find the average of the $I_g$ values in $supp(\phi_j^{l-1})$, denoted as $G_j$. If $G_j > \rho G_{mean}$, then we refine the particular B-spline basis according to the refinement procedure explained in section 2.2. The threshold value, $\rho$, is set according to the amount of refinement required to get the desired accuracy. A higher value of $\rho$ results in less control grid refinement.
3. After refinement, we collect all the active B-spline basis functions, control points and cells to get the hierarchical basis functions $\phi_{thb}^l$, control points $P_{thb}^l$ and cells $\Omega^l$ respectively. $\mathbf{M}^l$ matrix is computed for the refined control grid.
4. Begin the loop for iterations, $t = 1, 2, ..., i_{max}$, where $i_{max}$ is the maximum number of iterations for a particular level.
   (a) Compute the matrix consisting of the energy functional values, $\delta\mathbf{E}^{l,t}(f(\mathbf{x}))$.
   (b) We update the control points by solving $\mathbf{M}^l(\mathbf{P}_{thb}^{l,t+1} - \mathbf{P}_{thb}^{l,t}) = -\epsilon(\delta\mathbf{E}^{l,t}(f(\mathbf{x})))$.
   (c) Using the updated control points, the spatial transformation is computed and we have $f(\mathbf{x}) = \sum_{j=1}^{N_{thb}} P_{thb,j}^{l,t+1} \phi_{thb,j}^l$, where $N_{thb}$ is the total number of THB-spline basis functions.
   (d) Using the spatial transformation $f(\mathbf{x})$, we obtain image $I_1^l(f(\mathbf{x}))$.
5. The resulting image is used as the source image for the next refinement level.
6. The process terminates upon reaching the desired accuracy, that is, the best possible match between the source and target images.

In Figure 1, we demonstrate the application of image registration using adaptive local refinement on a synthetic image. We conduct the registration on five refinement levels with the number of control points of $529$, $1,159$, $2,611$, $5,935$ and $15,019$ respectively. The corresponding RS values are $55.18\%$, $92.17\%$, $94.83\%$, $97.41\%$ and $98.60\%$. From the results, we can observe that large deformations are captured on a coarser grid and finer deformations are captured on finer grids. Only regions with large change in deformation are refined. We thus make the entire registration process more efficient and fast, maintaining the same level of accuracy. To demonstrate the advantages of using THB-splines, we check the structure of the matrix $\mathbf{M}$ after five levels of refinement. As shown in Figure 2, the sparsity of the resulting matrices is improved by nearly $26\%$ compared to HB-spline basis functions. This improvement is significant in conducting the FEM-based implementation. By using THB-spline basis functions we introduce local refinement and fewer control points, hence reducing the size of matrices. In addition, the improvement of matrix sparsity reduces the computational cost and makes the entire registration process more efficient.

In summary, adaptive local refinement using THB-splines basis functions has the following advantages:

- We reduce the computational cost by introducing local refinement as compared with using uniformly refined grids. The desired level of accuracy is reached using fewer number of control points and at the same time we capture fine deformations on locally refined grids which cannot be captured properly on coarser grids.

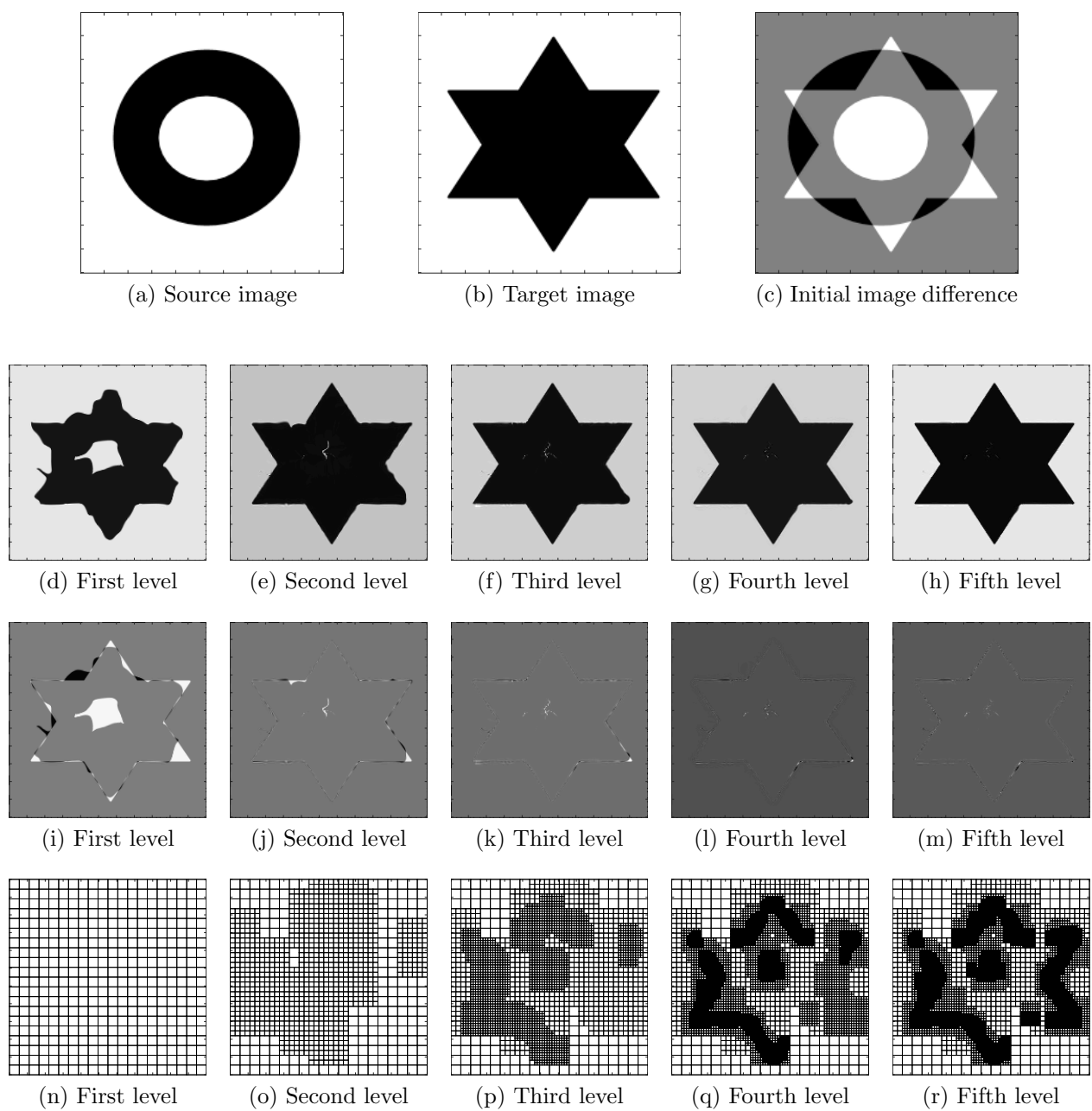

Figure 1: Registration of ring to star: The initial source image, the target image and initial differences are shown in (a-c). The evolving images after each stage of registration are shown in (d-h). The differences between the target and evolving images at each level are shown in (i-m). Meshes generated at different levels based on THB-spline basis refinement are shown in (n-r).

- THB-spline basis functions increase the sparsity of the matrices generated in FEM

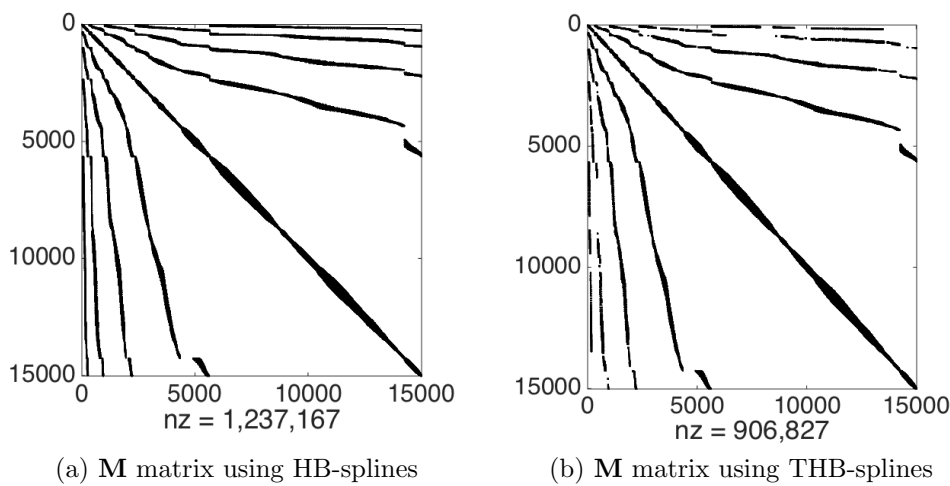

(a) **M** matrix using HB-splines  (b) **M** matrix using THB-splines

Figure 2: Comparison of the **M** matrix structure after five levels of refinement using HB-splines (a) and THB-splines (b) for the example shown in Figure 1. The sparsity ratio (THB-splines/HB-splines) is 73.30%.

based methods. As the computational cost depends on the number of non-zero entries of these matrices, this helps in making the solver more efficient.

## 4. NUMERICAL RESULTS AND DISCUSSIONS

The efficiency of the proposed registration framework is evaluated on two pairs of synthetic images and three pairs of medical images. The evolving images at the end of each refinement level are recorded. The difference between the evolving image and the target image at each refinement level is illustrated, along with the refined grids. All the results were generated on a computer with 2.5 GHz quad-core Intel Core i7 processor and 16GB RAM. Finally, a comparison of the similarity metrics and the computational time are shown in Tables 1 and 2 along with the results obtained from the uniform grid refinement technique.

We test the performance of our method on two pairs of synthetic images Figures 3-4 and three pairs of medical images Figures 5-7. All the numerical examples are implemented using the adaptive grid refinement strategy described in section 3. The regularization parameters $\lambda_1$ and $\lambda_2$ are set as 0.0001 for all our numerical examples. The time step $\epsilon$ is set by observing the convergence of the similarity ratio for the evolving images at each time step. At a particular level, when the values of $RS$ do not improve much or even start to decrease, the grid is refined locally. The control points are added in the areas of large change in deformation. The control mesh is reinitialized and the resulting image obtained at the previous refinement level is used as the source image for the newly refined grid. The registration process is not sensitive to the value of $\epsilon$ at a particular refinement level, but for different refinement levels the finer grids require $\epsilon$ to be comparatively smaller to accurately capture the finer deformations. The value of $\rho$ depends on the complexity of the image being studied. For medical images with more complex information, it is necessary to use more refinement levels and a larger value of $\rho$ for the finer refinement levels to prevent addition of too many control points. At the start of the registration process, these values are set for each refinement level, and then maintained fixed during the registration process.

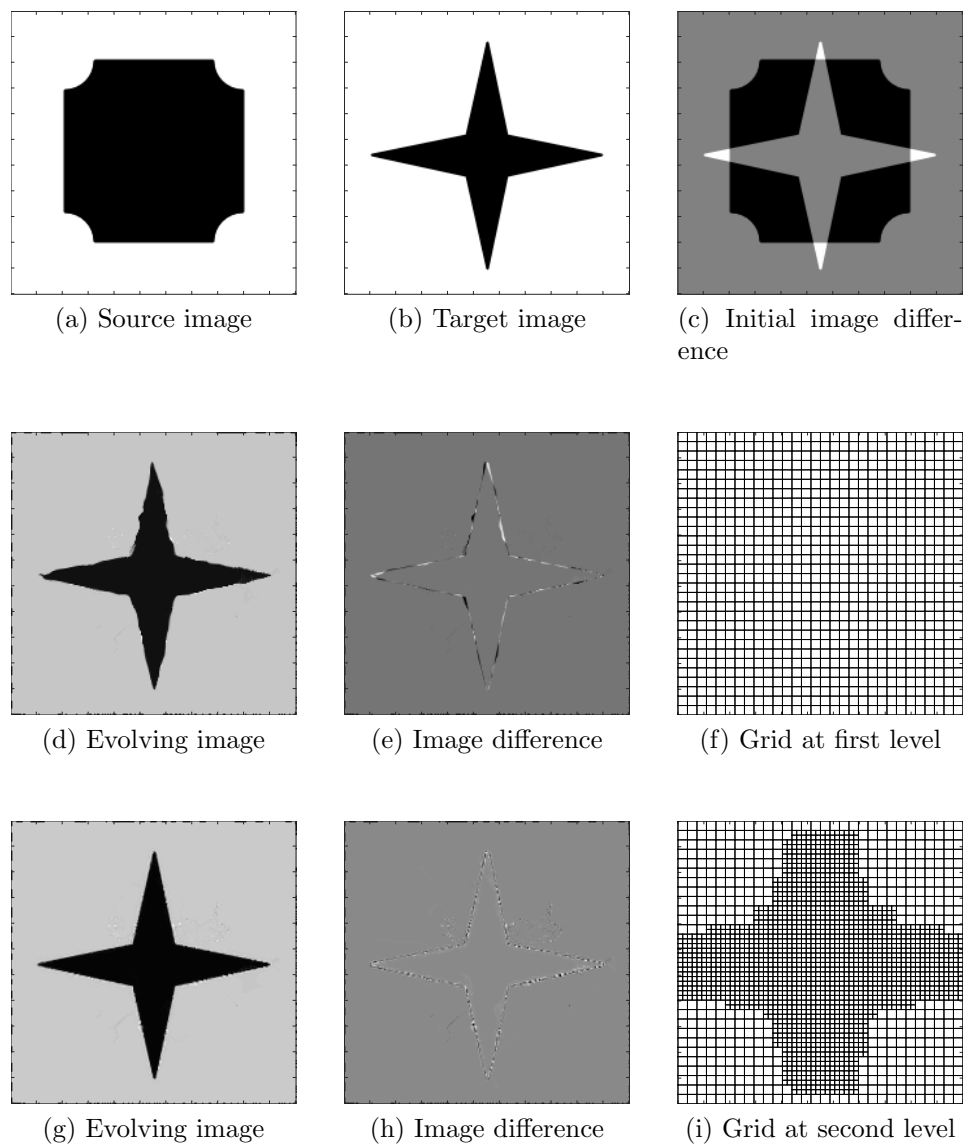

Figure 3: Square to Star: The initial source image, the target image and initial differences are shown in (a-c). The evolving images, differences between the target and evolving images and meshes generated at different levels based on THB-spline basis refinement for the first and second refinement levels are shown in (d-f) and (g-i) respectively.

For Figure 3, the values of $\epsilon$ are set as 0.025 and 0.005 for the first and second refinement levels respectively. The value of $\rho$ is set to 1.5. To achieve smooth convergence, we set $\epsilon$ to a lower value when conducting the registration on the finer levels. This is done to gradually deform the image and capture the highly localized deformations accurately. In Figure 4, we demonstrate the registration of two images with very different topology. To obtain the desired accuracy, we conduct the registration on three refinement levels. The values of $\epsilon$ and $\rho$ are set as 0.015 and 1.5 for the first two levels, 0.005 and 2 for the third level, respectively. For the finest level we set the value of $\rho$ to a higher value in order to avoid introducing too many control points in the image domain. Thus we achieve a lot of flexibility in terms of the amount of refinement required. Depending upon the nature of the registration problem, we can easily control the local refinement of the grid.

Table 1: Comparison of $MSD$ and $RS$ using our method (THBS) with uniform B-spline refinement (UBS).

| Image | Level | Iterations | | Control Points | | $MSD$ | | $RS$ | |
|---|---|---|---|---|---|---|---|---|---|
| | | UBS | THBS | UBS | THBS | UBS | THBS | UBS | THBS |
| Square-Star (Fig. 3) | 1 | 25 | 25 | 1,089 | 1,089 | 128.11 | 128.11 | 97.20 | 97.20 |
| Image Size: $220 \times 220$ | 2 | 6 | 6 | 3,969 | 2,382 | 20.68 | 23.03 | 99.15 | 99.12 |
| Sun-Four Stars (Fig. 4) | 1 | 15 | 15 | 1,089 | 1,089 | $1.70 \times 10^3$ | $1.70 \times 10^3$ | 69.64 | 69.64 |
| Image Size: $220 \times 220$ | 2 | 20 | 20 | 3,969 | 2,577 | 181.57 | 174.76 | 95.92 | 95.93 |
| | 3 | 7 | 10 | 15,129 | 6,498 | 18.70 | 21.77 | 98.99 | 98.91 |
| Brain 1 (Fig. 5) | 1 | 10 | 10 | 1,089 | 1,089 | 216.60 | 216.60 | 30.05 | 30.05 |
| Image Size: $256 \times 256$ | 2 | 20 | 20 | 3,969 | 2,142 | 84.00 | 90.97 | 73.79 | 73.32 |
| | 3 | 10 | 20 | 15,129 | 5,697 | 41.38 | 45.42 | 82.89 | 82.60 |
| | 4 | - | 10 | - | 6,807 | - | 30.58 | - | 84.35 |
| Brain 2 (Fig. 6) | 1 | 5 | 5 | 1,089 | 1,089 | 215.80 | 215.80 | 38.15 | 38.15 |
| Image Size: $512 \times 512$ | 2 | 15 | 15 | 3,969 | 2,322 | 70.79 | 72.90 | 73.28 | 72.98 |
| | 3 | 15 | 20 | 15,129 | 6,261 | 33.80 | 37.44 | 81.71 | 80.33 |
| | 4 | - | 10 | - | 10,077 | - | 24.89 | - | 83.03 |
| Brain 3 (Fig. 7) | 1 | 10 | 10 | 1,089 | 1,089 | 259.59 | 259.59 | 30.36 | 30.36 |
| Image Size: $256 \times 256$ | 2 | 20 | 20 | 3,969 | 2,945 | 96.15 | 100.32 | 67.55 | 66.66 |
| | 3 | 20 | 20 | 15,129 | 9,642 | 46.14 | 48.79 | 82.63 | 81.73 |
| | 4 | - | 10 | - | 10,344 | - | 31.16 | - | 86.60 |

In Figures 5-7, we implement the registration process on medical images, specifically to brain MRI images. Unlike synthetic images, medical images have more complex features. It is therefore crucial that we ensure the smooth convergence of the solution, so that the final images obtained are able to represent all the features correctly. The registration is conducted on four refinement levels. We set the values of $\epsilon$ as 0.008 for the first two levels and 0.005 for the last two levels. For Figures 5-6, the value of $\rho$ is set as 1.5 for the first two stages of refinement and 3 for the last refinement stage. For Figure 7 we set $\rho$ as 1 for the first two stages of refinement and 3 for the final stage of refinement. From the results we observe that an accuracy of around 98% and 83% is achieved for synthetic and medical images respectively. The final images are matched well with the corresponding target images. We can also observe that the solver can detect the areas of image differences at each level and locally refine these regions. This step ensures that we only add control points where image deformation is expected and thus capture these localized changes efficiently.

We also perform the registration of the examples using uniform grid refinement strategy. The comparison of the values of MSD, RS, number of control points and the running time after each refinement level is shown in Tables 1 and 2. The following observations can be

Table 2: Comparison of the total CPU time in seconds using our method (THBS) with uniform B-spline refinement (UBS).

| Image | Level | Iterations | | Control Points | | CPU Time (seconds) | |
|---|---|---|---|---|---|---|---|
| | | UBS | THBS | UBS | THBS | UBS | THBS |
| Square-Star (Fig. 3) | 1 | 25 | 25 | 1,089 | 1,089 | 42.32 | 40.02 |
| Image Size: $220 \times 220$ | 2 | 6 | 6 | 3,969 | 2,382 | 54.75 | 38.01 |
| | | | | | | total: 97.07 | total: 78.03 |
| Sun-Four Stars (Fig. 4) | 1 | 15 | 15 | 1,089 | 1,089 | 26.77 | 29.85 |
| Image Size: $220 \times 220$ | 2 | 20 | 20 | 3,969 | 2,577 | 159.86 | 93.65 |
| | 3 | 7 | 10 | 15,129 | 6,498 | 575.81 | 226.73 |
| | | | | | | total: 762.44 | total: 350.23 |
| Brain 1 (Fig. 5) | 1 | 10 | 10 | 1,089 | 1,089 | 18.60 | 24.73 |
| Image Size: $256 \times 256$ | 2 | 20 | 20 | 3,969 | 2,142 | 160.54 | 86.82 |
| | 3 | 10 | 20 | 15,129 | 5,697 | 713.20 | 305.30 |
| | 4 | - | 10 | - | 6,807 | - | 252.93 |
| | | | | | | total: 892.34 | total: 669.78 |
| Brain 2 (Fig. 6) | 1 | 5 | 5 | 1,089 | 1,089 | 11.97 | 34.12 |
| Image Size: $512 \times 512$ | 2 | 15 | 15 | 3,969 | 2,322 | 125.36 | 126.18 |
| | 3 | 15 | 20 | 15,129 | 6,261 | 1,310.50 | 416.96 |
| | 4 | - | 10 | - | 10,077 | - | 616.88 |
| | | | | | | total: 1,447.83 | total: 1,194.14 |
| Brain 3 (Fig. 7) | 1 | 10 | 10 | 1,089 | 1,089 | 18.09 | 22.40 |
| Image Size: $256 \times 256$ | 2 | 20 | 20 | 3,969 | 2,945 | 155.05 | 108.04 |
| | 3 | 20 | 20 | 15,129 | 9,642 | 1,259.65 | 506.79 |
| | 4 | - | 10 | - | 10,344 | - | 297.80 |
| | | | | | | total: 1,432.79 | total: 935.03 |

drawn from the results:

- The maximum deformation in the image is captured by the coarser grids. By observing the RS value after the initial few refinement levels, we can see that the maximum increase in similarity takes place on coarser levels. A larger time step is chosen to accelerate the deformation. On finer levels, the increase in RS is slower as compared to the coarser levels. Here highly localized deformations are captured and the evolving image gradually deforms to match the target image.

- Compared to uniform refinement, our THB-spline based method is more efficient for the same level of accuracy. Although we need slightly more time to initialize the THB-spline data structure, we use fewer control points and reach the desired accuracy faster. In this way, we are able to reduce the size of the matrix and make the solver much more efficient.

- Compared to using HB-splines, THB-splines improve the sparsity of the $\mathbf{M}$ matrix, further improving the computational efficiency.

- The computational time taken for the registration process depends mainly on the complexity of the images and to a lesser extent on the size of the images. This is because to accurately capture the features to be registered, we need more refinement levels and refine more regions at each stage than for images with less complex information.

## 5. CONCLUSION

In this paper, we propose an efficient way to solve FEM-based image registration problems using local refinement and THB-spline basis functions. The registration process is conducted on a series of grids capturing both large and finer deformations inside the image efficiently without introducing too many control points. This also helps in reducing the matrix size, leading to faster computation. Compared to HB-spline based method, THB-spline basis functions improve the sparsity of the matrices and satisfy partition of unity. The proposed method is tested on 2D synthetic and medical images. The improvement in the efficiency is shown by comparing the results of adaptive grid refinement with uniform grid refinement. For our future work, we plan to extend our scheme to 3D medical images and target real medical applications. In terms of computational cost for 3D image registration, we can predict that the proposed method will enhance the efficiency even more as compared to uniform refinement because regions of the image with uniform intensity (such as the background) can be represented with very few control points. The computational cost mainly depends on the complexity of the images. Even if the resolution of the 3D images is not the same in all the three directions, we can suitably choose the initial control grid to be coarser in the direction of lower resolution.

## ACKNOWLEDGEMENTS

The medical images were provided from (http://overcode.yak.net/15). The research at Carnegie Mellon University was supported in part by NSF CAREER Award OCI-1149591. The research at Bauhaus University Weimar was supported in part by the ITN-INSIST and ERC-COMBAT funded by the EU-FP7.

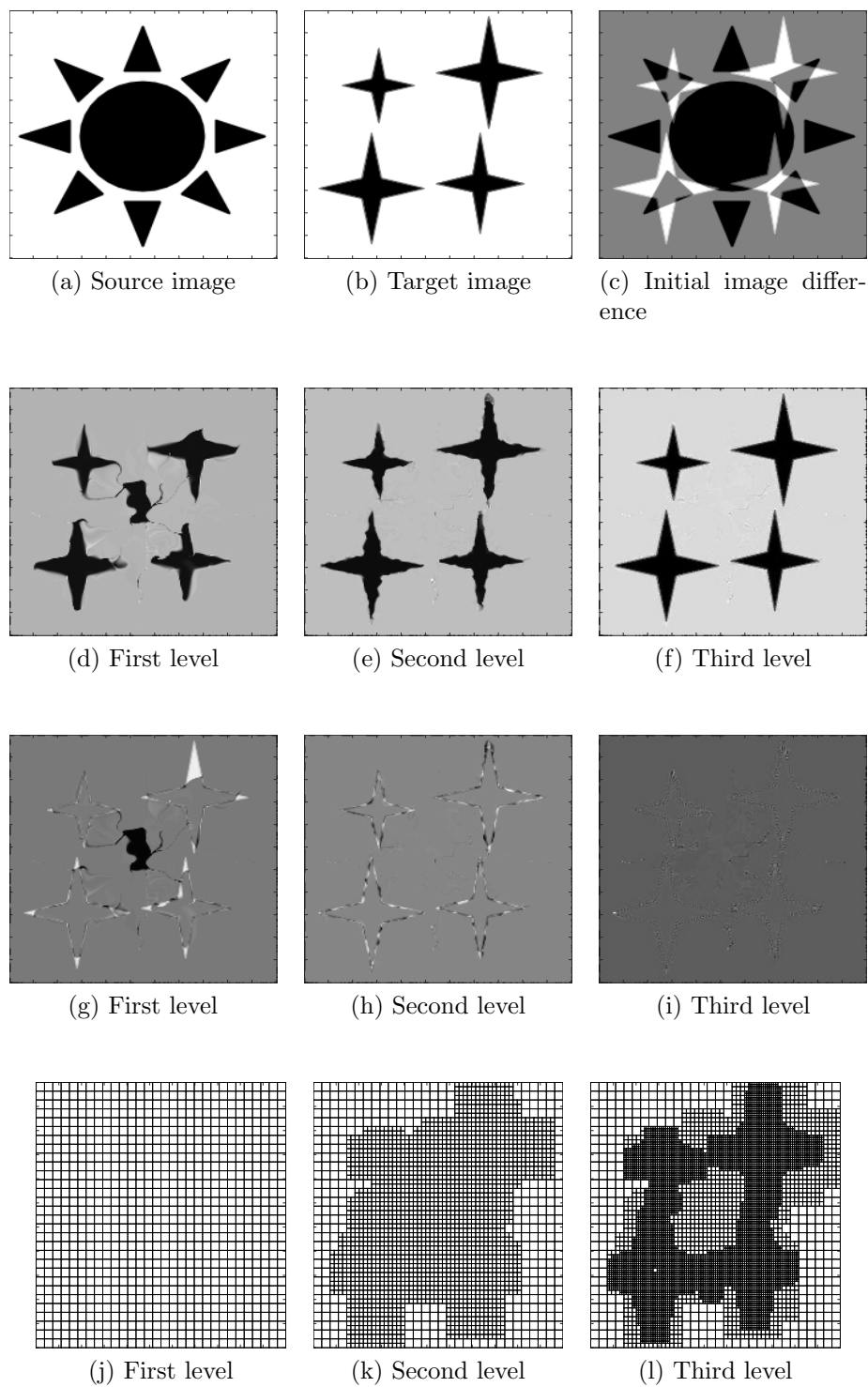

Figure 4: Sun to Four Stars: The initial source image, the target image and initial differences are shown in (a-c). The evolving images after each stage of registration are shown in (d-f). The differences between the target and evolving images at each level are shown in (g-i). Meshes generated at different levels based on THB-spline basis refinement are shown in (j-l).

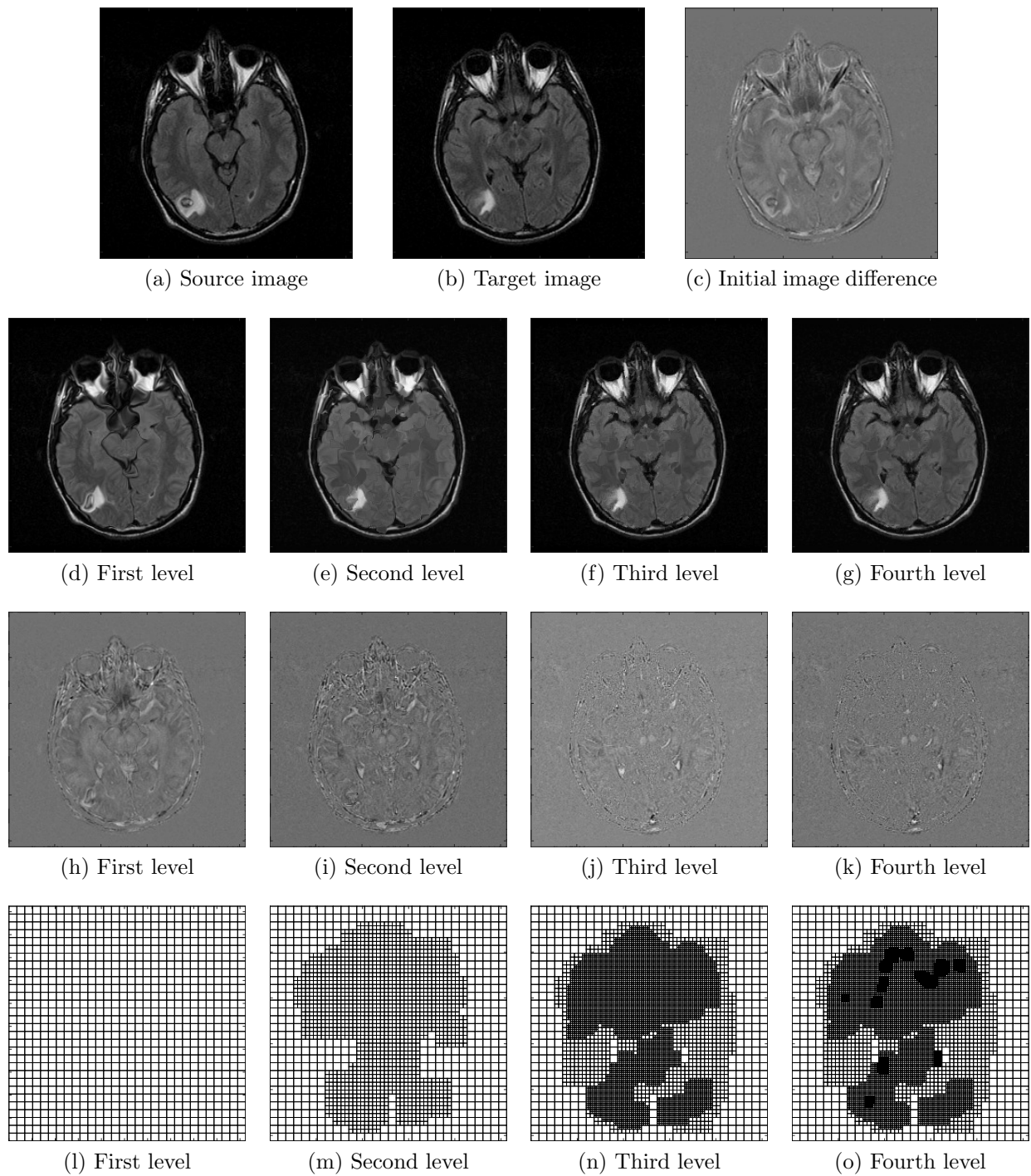

Figure 5: Brain 1: The initial source image, the target image and initial differences are shown in (a-c). The evolving images after each stage of registration are shown in (d-g). The differences between the target and evolving images at each level are shown in (h-k). Meshes generated at different levels based on THB-spline basis refinement are shown in (l-o).

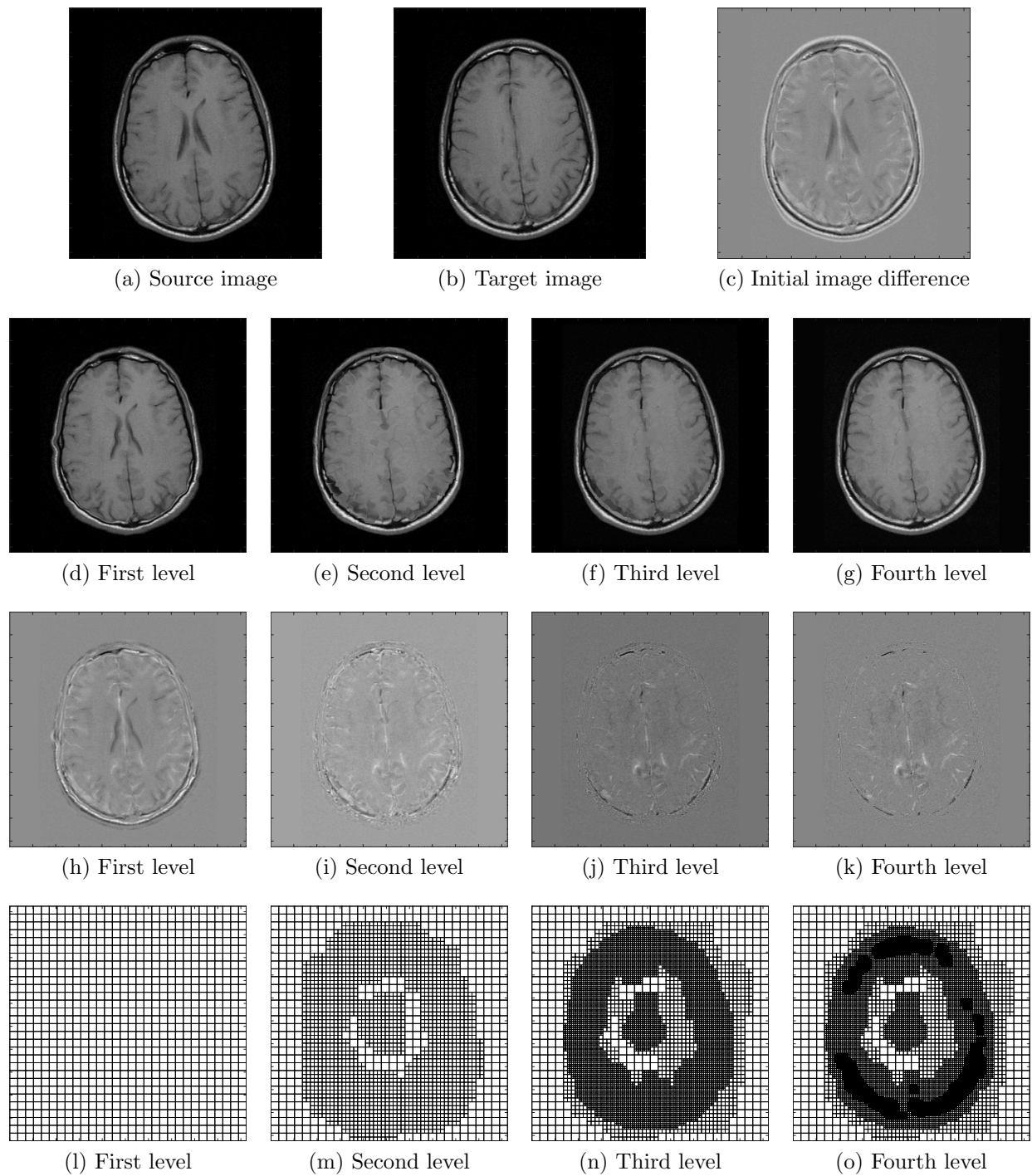

Figure 6: Brain 2: The initial source image, the target image and initial differences are shown in (a-c). The evolving images after each stage of registration are shown in (d-g). The differences between the target and evolving images at each level are shown in (h-k). Meshes generated at different levels based on THB-spline basis refinement are shown in (l-o).

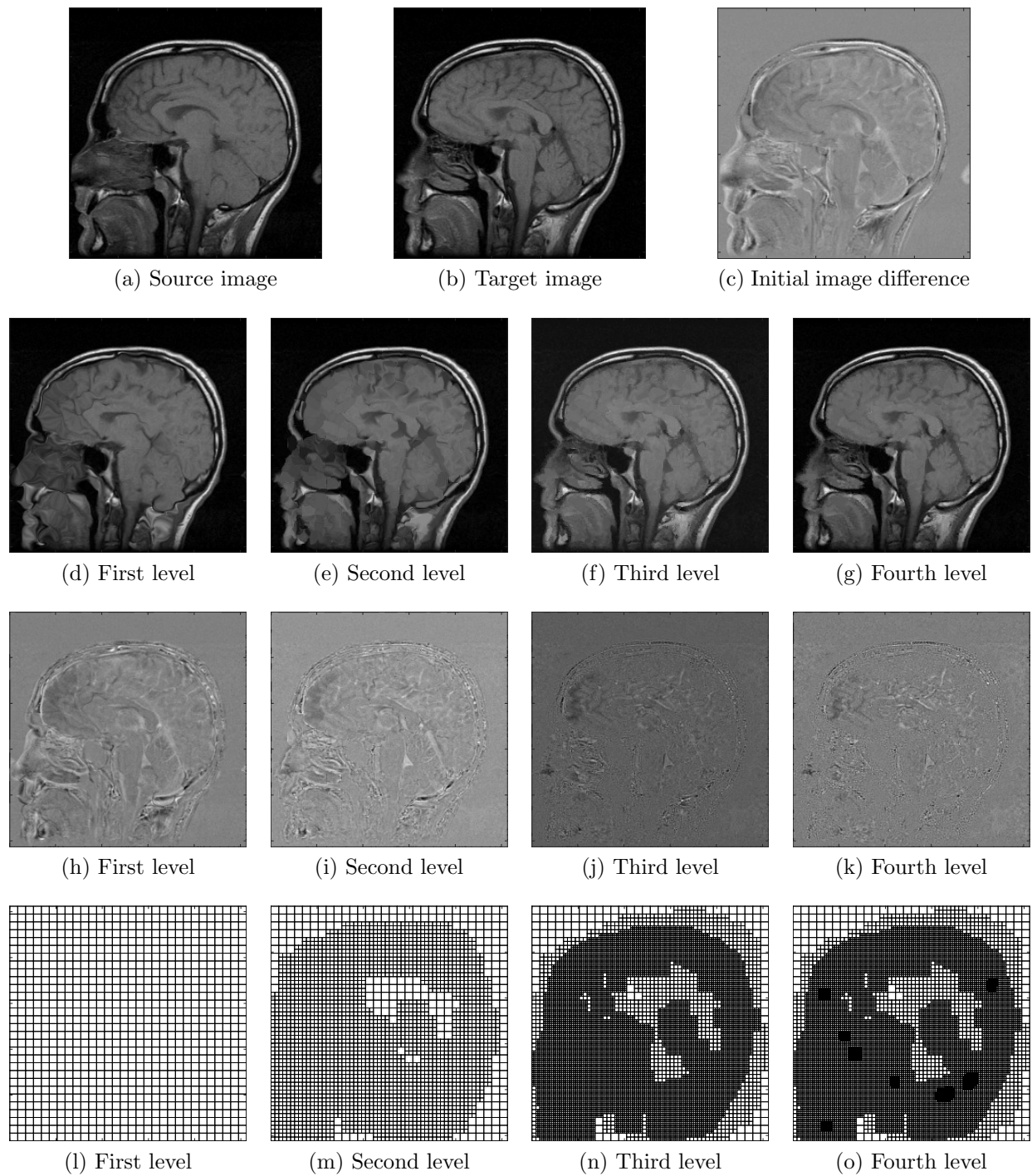

Figure 7: Brain 3: The initial source image, the target image and initial differences are shown in (a-c). The evolving images after each stage of registration are shown in (d-g). The differences between the target and evolving images at each level are shown in (h-k). Meshes generated at different levels based on THB-spline basis refinement are shown in (l-o).

19